\documentclass[12pt]{article}
\usepackage{amsmath,pifont}
\usepackage{amssymb}
\usepackage{amsthm}
\usepackage{psfig}
\usepackage[arrow,matrix,tips,2cell]{xy}
\usepackage[dvips]{graphics}
\newcommand{\ignore}[1]{\relax}

\newcommand{\C}{\mathbb C}
\newcommand{\R}{\mathbb R}
\newcommand{\Z}{\mathbb Z}

\newcommand{\re}{\operatorname{Re}}
\newcommand{\im}{\operatorname{Im}}

\newcommand{\Arg}{\operatorname{Arg}}

\newcommand{\ver}{\operatorname{Vert}}

\newcommand{\Alga}{\operatorname{Alga}}

\newtheorem{thm}{Theorem}

\newtheorem{lem}{Lemma}[section]
\newtheorem{lemmaA}{Lemma A\hspace{-5pt}}
\newtheorem{theoremA}{Theorem A\hspace{-5pt}}

\newtheorem{cor}{Corollary}

\newtheorem{coro}[lem]{Corollary}
\newtheorem{prop}[lem]{Proposition}

\theoremstyle{definition}
\newtheorem{defn}[lem]{Definition}

\newtheorem{exa}[lem]{Example}

\theoremstyle{remark}
\newtheorem{rmk}[lem]{Remark}

\newcommand{\tor}{(\C^\times)^{2}}

\newcommand{\rtor}{(\R^\times)^{2}}

\newcommand{\conj}{\operatorname{conj}}

\newcommand{\dd}{\partial}

\newcommand{\cp}{{\mathbb C}{\mathbb P}}
\newcommand{\rp}{{\mathbb R}{\mathbb P}}

\newcommand{\Log}{\operatorname{Log}}

\newcommand{\Vol}{\operatorname{Vol}}

\newcommand{\Int}{\operatorname{Int}}

\renewcommand{\setminus}{\smallsetminus}

\newcommand{\Area}{\operatorname{Area}}

\newcommand{\cP}{\mathcal P}
\newcommand{\cT}{\mathcal T}

\newcommand{\cS}{\mathcal S}
\newcommand{\cR}{\mathcal R}

\newcommand{\barE}{\overline{E}}
\newcommand{\barF}{\overline{F}}

\newcommand{\Pt}{\widetilde{P}}
\newcommand{\Rh}{\widehat{R}}
\newcommand{\Qt}{\widetilde{Q}}
\newcommand{\Rt}{\widetilde{R}}

\newcommand{\cb}{\textup{\ding{114}}}

\newcommand{\logf}{\div}

\newcommand{\apar}{\rightleftarrows}
\newcommand{\pline}{\mathbb{P}^1}
\newcommand{\bbP}{\mathbb{P}}

\DeclareMathOperator{\lcm}{lcm}
\DeclareMathOperator{\br}{br}

\begin{document}
\title
{Geometry of planar log-fronts}
\author{Grigory Mikhalkin and Andrei Okounkov}
\maketitle
\section{Introduction}

\subsection{Frozen boundaries}\label{motiv}

Given two polynomials $P(z,w)$ and $Q(z,w)$, one
can study solutions of the system
\begin{equation}
  \label{PQeq}
  \begin{cases}
     P(z,w)=0\,, \\
  Q(e^{-x} z, e^{-y} w) =0\,.
  \end{cases}
 \end{equation}
as a function of a point $(x,y)\in \R^2$.  A solution
$z(x,y)$ and $w(x,y)$
of \eqref{PQeq} solves a first order quasilinear PDE
\begin{equation*}
  \label{Burg}
  \frac{z_x}{z} + \frac{w_y}{w} = 1\,, \quad P(z,w)=0\,,
\end{equation*}
which is closely related to the complex Burgers equation
and arises in the theory of random surfaces, see \cite{Burg}.
The singularities of $z(x,y)$ and $w(x,y)$ occur when
\eqref{PQeq} has a multiple root, that is, when the two
curves in \eqref{PQeq} are tangent.
In the random surface context,
this marks the boundary between order (e.g.\ crystalline
facet) and disorder,
called \textit{frozen boundary}.

\begin{figure}[!h]
  \centering
  \scalebox{0.4}{\includegraphics{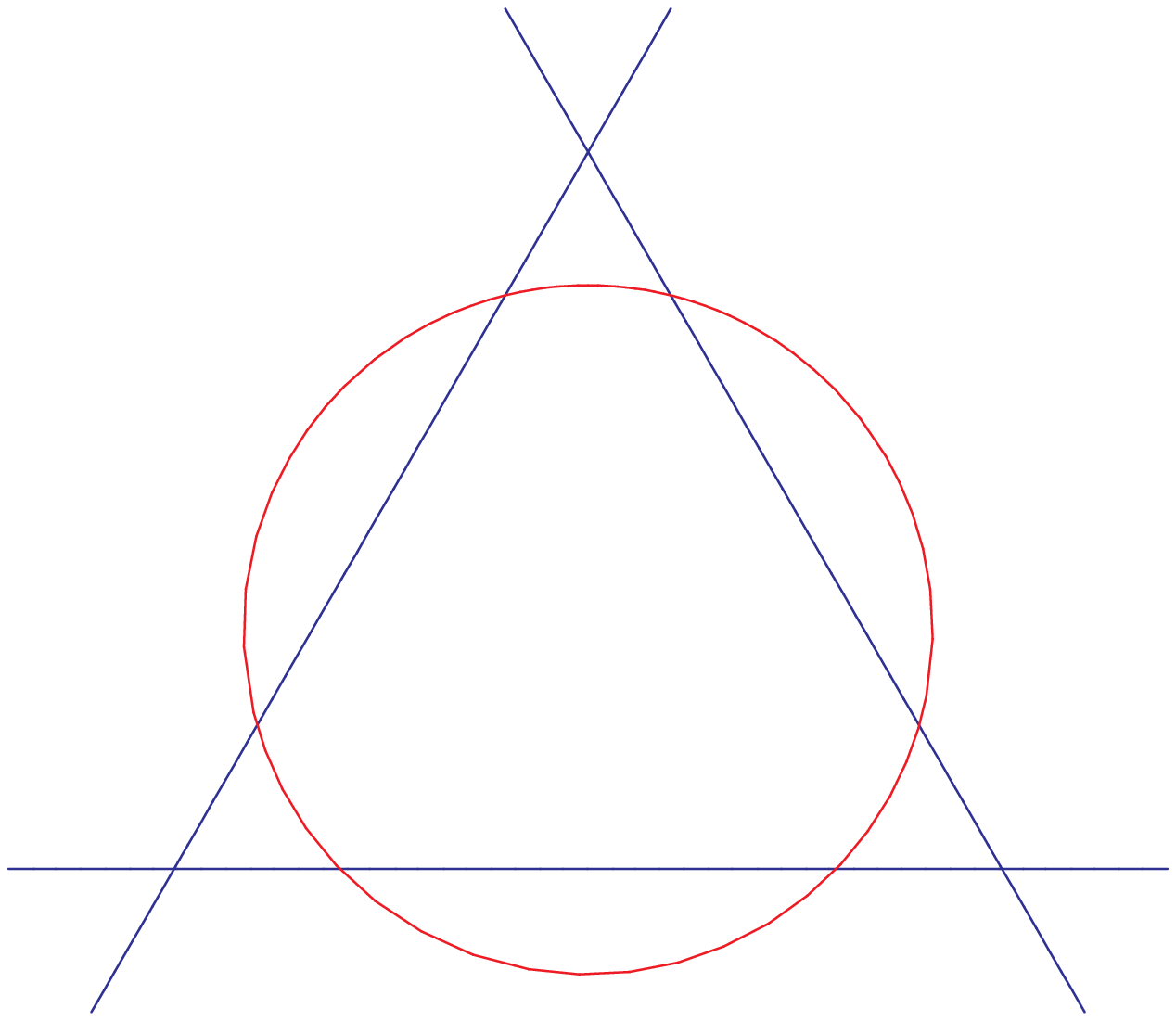}}
\quad \scalebox{0.4}{\includegraphics{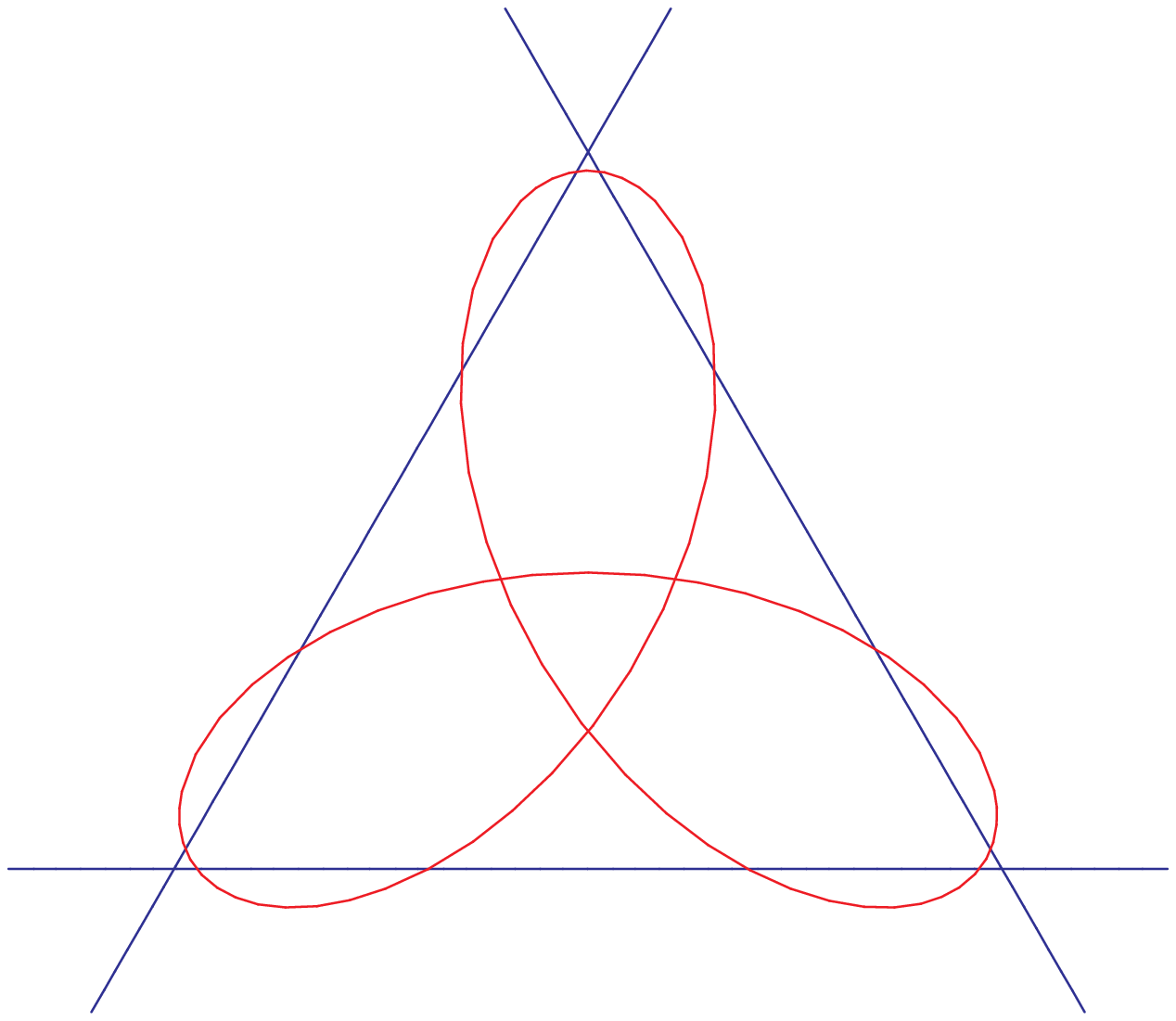}}
  \caption{Curves $P$ and $Q$ for the frozen boundary in Figure \ref{froz}}
  \label{PandQ}
\end{figure}

\begin{figure}[!h]
  \centering
  \rotatebox{90}{\scalebox{0.5}{\includegraphics{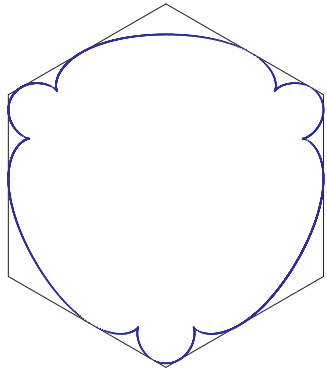}}}
  \caption{An example of frozen boundary}
  \label{froz}
\end{figure}

For example, let $P$ and $Q$ define
rational curves in $\rp^2$ of degree $2$ and $4$, respectively,
positioned with respect to the coordinate axes of $\rp^2$ as illustrated in
Figure \ref{PandQ}.
Points where $P$ and $Q$ intersect the axes may be
fixed so that the frozen boundary will be inscribed in a hexagon as
illustrated in Figure \ref{froz}.
The probabilistic meaning of this curve
is the following.

The well-known
arctic circle theorem of Cohn, Larsen, and Propp \cite{CLP}
states that the limit shape of a typical 3D partition contained
in a cube has a frozen boundary which is a
circle. More precisely, projected in
the $(1,1,1)$-direction, the cube becomes a hexagon and
the $(1,1,1)$-projection of the frozen boundary is an inscribed
circle. Now suppose we additionally weight each 3D partition $\pi$ by a weight
which a product over all boxes
$$
\cb=(i,j,k)\in \Z^3
$$
in $\pi$ of some periodic
function of $i-j$ and $j-k$. If the period equals $2$, the
frozen boundary will look as in Figure \ref{froz}. The coefficients
of the curve $P$ are obtained from the periodic weights, while
the coefficients of $Q$ are fixed by boundary conditions.
See \cite{Burg} for a detailed discussion of this procedure in general.

The geometry and, especially, the
singularities of frozen boundaries are of considerable
interest. Note that these are curves of some complexity:
in our current example frozen boundary has degree $28$
and genus $27$.

\subsection{Log-front}

This motivates the following definition.
Given two curves $P,Q\subset\tor$, consider
\begin{equation} \label{defR}
R =\left\{(a,b) \, \big| \,
\textup{$\tau_{a,b}\cdot Q$ is tangent to $P$}
\right\} \subset \tor
\end{equation}
where
$$
\tau_{a,b}: \tor\to\tor
$$
it the dilation $(z,w)\mapsto (az,bw)$ and
tangency at a singular point  means that for
some branches of the two curves in question
their unique tangent lines coincide.
For simplicity, we assume that
tangency occurs only at isolated points, that is,
no component of $P$ is a dilate of a component of $Q$.

We call
the curve $R$ the \emph{log-front} of $P$ and $Q$.
The reason for such name will be explained below.
We will denote the construction \eqref{defR}
by
$$
R=P\logf Q\,.
$$
The obvious property
$$
Q\logf P = R(a^{-1},b^{-1})
$$
and \eqref{double_dual} below justify such notation.

Clearly, $R$ is an algebraic curve.
Various ways to compute its equation
$R(a,b)=0$ will be discussed in Section \ref{sRhat}.
The frozen boundary from Section \ref{motiv} is
given by $R(e^x,e^y)=0$.

\subsection{Symmetry between $Q$ and $R$}
The definition of $R$ can be recast into several
equivalent forms.
Let $\cP$ denote the hypersurface
$$
\cP=\left\{P(az,bw)=0\right\} \in \tor \times \tor
$$
and let $\pi_{ab}$ and $\pi_{zw}$ be the projections
from $\tor\times \tor$ to the respective factors.

The tangency in definition \eqref{defR} can be
rephrased by saying that $R$ is formed by critical
values of the map
\begin{equation}
  \label{cS}
    \pi_{ab}: \quad\cS=\cP \cap \pi_{xy}^{-1} \left(Q\right) \to \tor
\end{equation}
of complex surfaces.
In other words, the surface
$$
\pi_{ab}^{-1} \left(R\right) \cap \pi_{xy}^{-1} \left(Q\right)
\cong R \times Q
$$
is tangent to the hypersurface $\cP$ along their intersection, showing
a certain symmetry
between the roles of $Q$ and $R$ for fixed $P$. In
particular
\begin{equation}
  \label{double_dual}
 Q \subset P \logf(P \logf Q)\,.
\end{equation}
The multiplicity, with which $Q$ occurs in the
right-hand side of \eqref{double_dual} equals
the degree of the
\emph{logarithmic Gau\ss\ map} of $P$, see below.

Yet another way to say the same thing is
that $R$ is the envelope of the family
of curves
\begin{equation}
  \label{envelope}
     \big\{P(az,bw)=0\big\}_{Q(z,w)=0}\,,
\end{equation}
indexed by points of $Q$

\subsection{Classical constructions}

Among examples of the operation \eqref{defR}
 there are the
following two classical constructions.

First, let $P$ be a general line, for example,
$$
P(z,w) = z + w + 1 \,.
$$
Then its dilates $\tau_{ab}\cdot P$ form an open
set of the dual projective plane and hence $R$ is
an open set of the dual curve $Q^\vee$, namely, its
intersection with $\tor$. In this
case, \eqref{double_dual} becomes an equality.

Second, the additive analog of \eqref{defR}
\begin{equation} \label{defadd}
P - Q =\left\{(a,b) \, \big| \,
\textup{$Q(z-a,w-b)=0$ is tangent to $P$}
\right\} \subset \C^2
\end{equation}
is a limit case of \eqref{defR}. For
analytic plane curves, constructions $P\logf Q$ and \eqref{defadd}
are, in fact, equivalent by taking the logarithms.
When $P$ is
a circle of radius $t$
$$
P(z,w) = z^2 + w^2 - t^2
$$
the real locus of $P-Q$ contains
%by \eqref{envelope}
the front at time $t$ of a wave that was emitted at time zero from all
real points of $Q$ and is propagating with unit
velocity. Such curve is called a \emph{wave-front}.
It is this example that motivates
the general term \emph{log-front}.

Pl\"ucker formulas relate singularities of a curve $Q$
to the singularities of its dual $Q^\vee$. A formula
of Klein further constraints the singularities of the
real locus of $Q^\vee$. Analogous formulas for wave
fronts were obtained by O.~Viro in \cite{V-usp}. The
goal of this note is prove an analog of Pl\"ucker
and Klein formulas in the general case.

In probabilistic applications, the curve $P$ is
a real algebraic curve of a very special kind, namely,
it is a \emph{Harnack curves}. See \cite{Mi-obzor} and the Appendix
for a discussion of the properties of Harnack curves
and \cite{KOS,KO1} for connections with probability. Harnack
curves have many remarkable features that general real
plane curves lack. As it turns out, the assumption that
$P$ is Harnack is also essential for our
derivation of Klein-type formulas
for $R$.

It may
be noted here that both our formulas and their
proofs involve nothing but elementary geometry of
plane curves and would have been, no doubt,
obtained by Pl\"ucker and Klein had their seen
a need for them.

\subsection{Acknowledgments}

We are grateful to C.~Faber, R.~Kenyon and O.~Viro
for many useful discussions. Some of our formulas generalize
unpublished results obtained jointly
with C.~Faber.%, see \cite{FO}.

\section{Preliminaries}

There are many good books on plane algebraic
curves, see e.g. \cite{B,W}. With the probability audience
in mind, we collected in this section an explanation
of some basic
notions that will be used later.

\subsection{Newton polygons}\label{torsurf}

Let $P\subset\tor$ be a curve with the equation
$$
P(z,w)=\sum\limits_{i,j}p_{ij}\, z^i \, w^j=0\,.
$$
By definition, its \emph{Newton polygon} is
$$
\Delta_P=\operatorname{Convex\ Hull}
\left\{(i,j)\in\Z^2 \,\big|\, p_{ij}\neq 0\right\}\,.
$$
This is a refinement of the degree of $P$.
The polygon $\Delta_P$ is defined only  up to
translation by a lattice vector once we treat $P\subset\tor$
as a geometric curve rather than a polynomial. Let
$|\Delta_P|$ denote the number of lattice points in
$\Delta_P$.

The Newton polygon  defines a projective toric surface
\begin{equation}\label{cT}
\cT=\overline{\Big\{\big(z^i w^j\big)_{(i,j)\in \Delta_P}\Big\}}_{(z,w)\in \tor}
\subset \cp^{|\Delta_P|-1}
\end{equation}
of degree
$$
\deg \cT = 2 \Area(\Delta_P)\,.
$$
{}From now on, $P$ will denote the closed curve in
$\cT$ defined by $P(z,w)=0$. This is a
hyperplane section avoiding the torus
fixed points, which are the only possible singularities
of $\cT$. We call the components $\mathcal{D}_E$ of
$$
\partial \cT = \cT \setminus \tor = \bigcup_{\textup{edges $E$ of $\Delta_P$}} \mathcal{D}_E
$$
the \emph{boundary divisors} and the points of
$$
\partial P = P \cap \partial \cT
$$
the \emph{boundary points} of $P$. For simplicity
we assume throughout the paper that the boundary
points of both $P$ and $Q$ are smooth.

By construction,
$P$ intersects $\mathcal{D}_E$  in $|E|$ points
counting multiplicity, where $|E|$ is the length
of $E$ in lattice units. The multiplicities of these
intersection points define a partition $\lambda(E)$
of $|E|$ for every edge $E$ of $\Delta_P$. We will
call a lattice polygon with such additional
partition data a \emph{marked polygon}.

\subsection{Amoebas}

The image $\Log(P)$ of a curve $P$ under the map
$$
\Log((z,w)) = (\log|z|,\log|w|)
$$
is called the \emph{amoeba} of $P$. The boundary
points of $P$ give rise to the so-called
\emph{tentacles} of the amoeba $\Log(P)$, see
for example Figure \ref{fig_tang_amoeb} which shows
logarithmic images $\Log(\R P)$ of the real
loci of two plane curves.

\begin{figure}[!htbp]
  \centering
  \scalebox{0.33}{\includegraphics{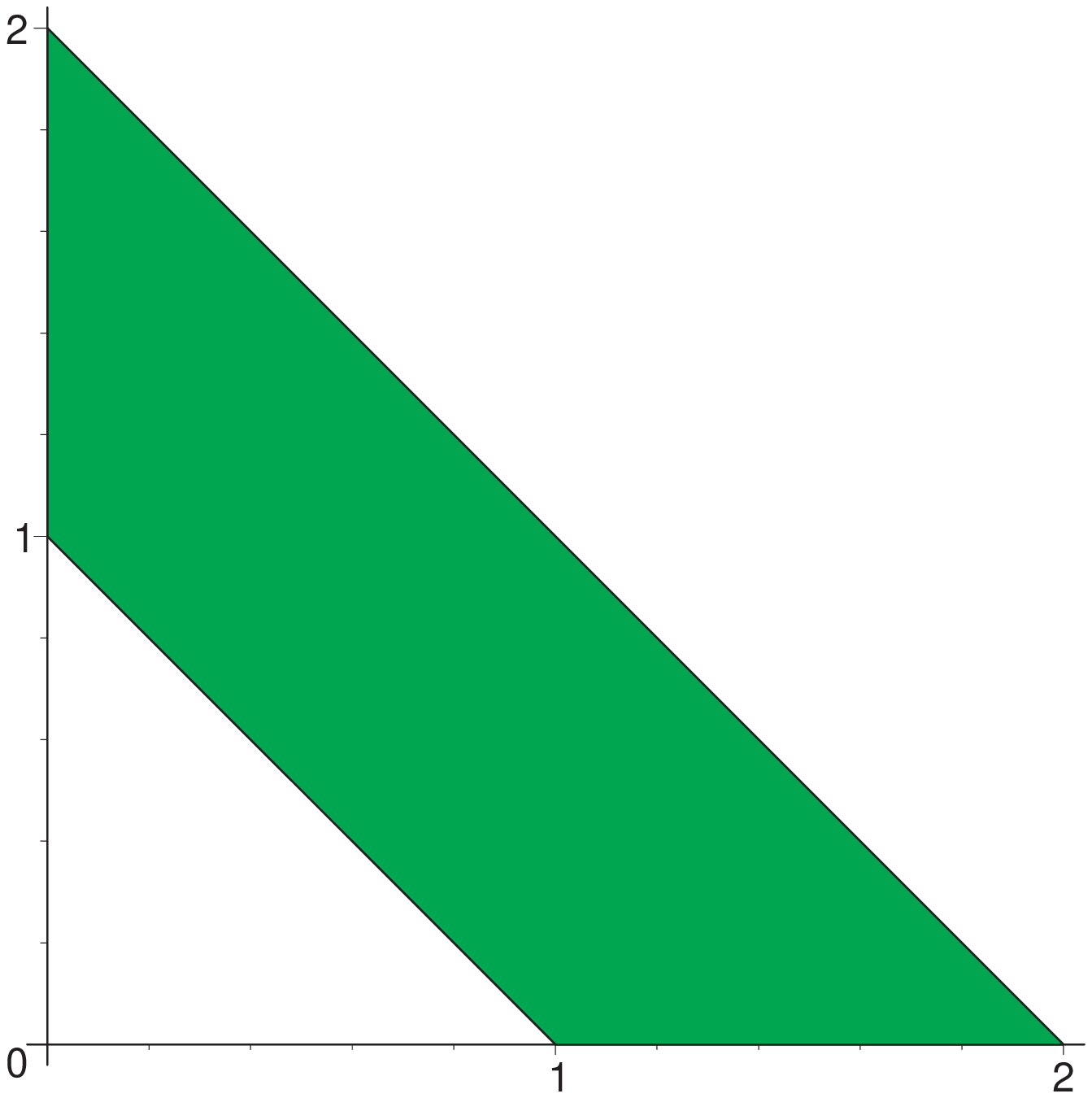}}
  \quad
  \scalebox{0.33}{\includegraphics{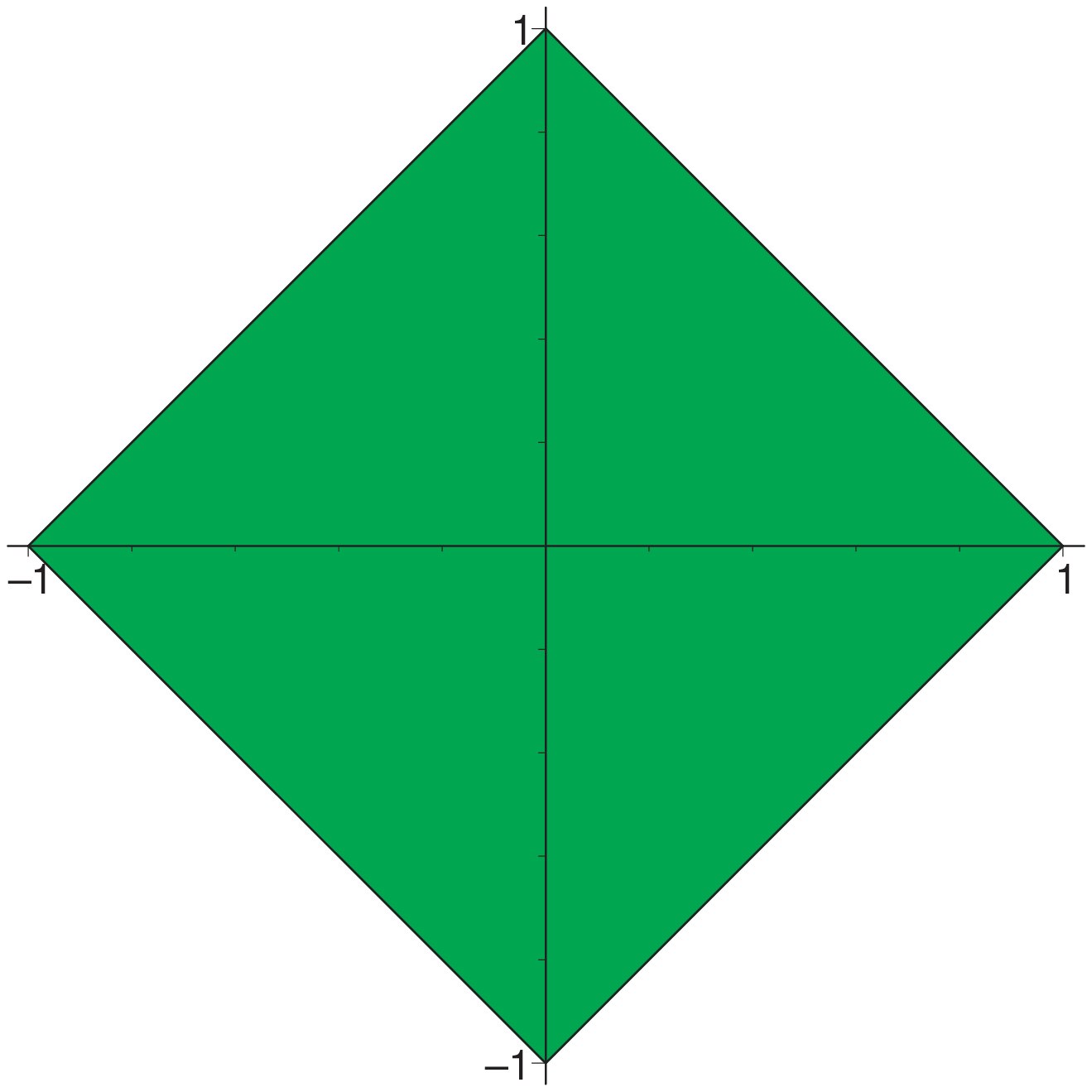}}
  \caption{Newtons polygons of curves in Figure \ref{fig_tang_amoeb}}
  \label{fig_Polygons}
\end{figure}

The directions of
the tentacles are the outward normals to the
sides of the Newton polygon. For example,
Newton polygons of curves from
Figure \ref{fig_tang_amoeb} are plotted in
Figure \ref{fig_Polygons}.
The number of
tentacles in a given direction is the number
of points of $P$ on the corresponding boundary
divisor.

\subsection{Logarithmic Gau\ss\ map}\label{subsectionGauss}

Since $\tor$ is an abelian group, the tangent spaces at
all points of $\tor$ are canonically identified with
tangent space $T_{(1,1)}$ at the identity $(1,1)\in\tor$. Mapping
the tangents to $P$ to their images in $\bbP(T_{(1,1)})$
gives the \emph{logarithmic Gau\ss\ map}
$$
\gamma_P: \Pt \to \pline\,,
$$
where $\Pt$ is the normalization of $P$. In coordinates,
$$
\gamma_P((z,w)) = \frac{z \frac{\partial}{\partial z} P}
{w \frac{\partial}{\partial w} P}\,.
$$
Note that the values of $\gamma_P$ at the boundary points of $P$
are determined by the slopes of the corresponding edges of $\Delta_P$.

The degree of  $\gamma_P$ can be computed as follows,
cf.\ \cite{Ka,Mi-obzor}.
Recall that the \emph{multiplicity} $m(p)$ of a point $p\in P$ is
the multiplicity with which $P$ intersects a generic
line through $p$. Another important characteristic of
a singular point is its \emph{Milnor number} $\mu(p)$, see \cite{Milnor}.
It may be defined as
the local intersection number of $\frac{\partial}{\partial z} P=0$
and $\frac{\partial}{\partial w} P=0$ at $p$.

\begin{prop}\label{dGauss}
We have
\begin{equation}
  \label{deg_gam}
     \deg \gamma_P = 2 |\Delta^\circ_P | + |\partial P| - 2 -
\sum_{p\in P}(\mu(p)+m(p)-1)\,.
\end{equation}
where $|\Delta^\circ_P |$ is the number of points in the
interior of the Newton polygon $\Delta_P$ and  $|\partial P|$
is the number of boundary points of $P$ not counting
multiplicity.
\end{prop}

Clearly, only
singular points of $P$ contribute to the sum over $p\in P$
in \eqref{deg_gam} since for a nonsingular point
$p$ we have $m(p)=1$ and $\mu(p)=0$.

\begin{proof}
If $P$ is smooth and transverse to
the boundary
then the proposition follows from
Kouchnirenko's formula \cite{Ku} (also using
Pick's formula for the area of a lattice polygon).

A singular point of $P$ subtracts $\mu(p)+m(p)-1$,
which is
the same computation as for the ordinary Gau\ss\ map,
see e.g.\ Theorem 7.2.2.\ in \cite{W}. The
effect of tangency to the boundary is determined similarly,
see e.g.\ Section \ref{proofPolygon} below.
\end{proof}

The critical points of $\gamma_P$ are known as the
logarithmic inflection points. Note that, for example,
a cusp is not a logarithmic inflection point, but a
sharp, or ramphoid, cusp (locally looking like $x^2=y^5$) is.
By the Riemann-Hurwitz
formula applied to $\gamma_P$, the count of
logarithmic inflection points with multiplicity is
\begin{equation}\label{logFlex}
\textup{\# logarithmic inflection points}=2\deg \gamma_P - \chi(\Pt).
\end{equation}

\subsection{Geometric genus and adjunction}

The geometric genus of $P$ is defined as the genus
$g(\Pt)$ of its normalization $\Pt$. It may
be computed as follows, cf. \cite{Kh},
\begin{equation}
  \label{adju}
  g(\Pt) = |\Delta^\circ_P | - \frac12 \sum_{p\in P} (\mu(p)+ \beta(p) -1)
\end{equation}
where $\beta(p)$ is the number of branches of $P$ through $p$.
%Thus we have
This is known as the adjunction formula.

The number
$\frac12 (\mu(p)+ \beta(p) -1)$ is a nonnegative integer which
vanishes if $p$ is a smooth point of $P$. It equals $1$
for cusps and nodes, so if the curve has only those
singularities, we get
\begin{equation}
\textup{\# cusps} +
\textup{\# nodes}= |\Delta^\circ_P | - g(\Pt) \label{n_nodes}\,.
\end{equation}
Combining \eqref{deg_gam} and \eqref{adju} we get
\begin{equation}
\sum_{p\in P} (m(p)-\beta(p))=
- \deg \gamma_P - \chi(\Pt\setminus \partial P)\label{n_cusps}\,.
\end{equation}

\begin{rmk}
A geometrically-minded reader might appreciate the following
alternative proof (or rather a rephrasing of the proof)
of \ref{n_cusps}.

Suppose that $\Pt\setminus \partial P$ is smooth.
We may deduce the equation
\begin{equation}\label{eq-smooth}
\deg \gamma_P=- \chi(\Pt\setminus \partial P)
\end{equation}
by the following application of the maximum principle, cf. \cite{Mi}.
Let $\lambda:\R^2\to\R$ be a linear map which we may choose
with a generic slope of the kernel.
The function $\lambda\circ\Log$ is a pluriharmonic function on $\tor$
and thus restricts to a harmonic function $h$
on $\deg \gamma_P=- \chi(\Pt\setminus \partial P)$.
Thus all critical points of $h$ are of index 1.

Furthermore, $h$ exhibits $P$ as a cobordism between $P_-\subset\dd P$
and $P_+\subset\dd P$, where the splitting $\dd P=P_+\cup P_-$
corresponds to $h^{-1}(+\infty)$ and $h^{-1}(-\infty)$.
Thus, the number of critical points of $h$ equals $\chi(P-P_+,P_-)=\chi(P\setminus\dd P)$
(the latter equality follows from additivity of Euler characteristic).
On the other hand the number equals $\deg\gamma_P$, so we get \eqref{eq-smooth}.
Each singularity $p$ of $P_-\subset\dd P$ subtracts $\sum(\mu(P)-\beta(p)-1)$
from $-\chi(P)$ and $\sum(\mu(P)-m(p)-1)$ from $\deg\gamma_P$.
\end{rmk}

\subsection{Nodal and cuspidal numbers of a plane curve}\label{bcdefn}
\begin{figure}[!htbp]
  \centering
  \scalebox{0.4}{\includegraphics{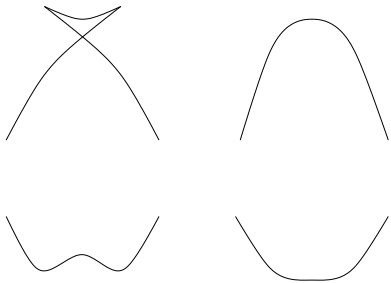}}
  \caption{Two cusps coalescing in a $E_6$ singularity.}
  \label{E6}
\end{figure}

The right-hand side of \eqref{n_cusps}
may be interpreted as the total
\emph{cuspidal number} $c(P)$ of $P$. It equals the
number of cusps if $P$ has only cusps and
ordinary multiple points. An example of
formula \eqref{n_cusps} may be seen in Figure
\ref{E6}. In the bottom half of Figure \ref{E6},
we see two inflection point come together to
form a point of the form $y=x^4$. For the
dual curves, which have the same genus and the same
degree of the
logarithmic Gau\ss\ map, we see a merger of
cusps into a $E_6$ singularity $y^3=x^4$ which
has the cuspidal number equal to 2.

Another useful number is the {\em nodal number} $b(P)$.
%Let $s\in \bar{A}$ be a singular point and $B_1,\dots,B_l$
%be the branches of $\bar{A}$ that pass through $s$.
%Denote with $\#(B_j,B_k)$ the local intersection number
%of the branches $B_j$ and $B_k$ at $s$.
First we define the {\em local nodal number}
{$b_s(P)$} at every singular point $s\in P$
as the sum $\sum\limits_{j\neq k}\#(B_j,B_k)$ of the intersection
numbers over all distinct pairs of branches $B_j$, $B_k$ via $s$.
Clearly, $b_s(P)=0$ if and only if $s$ is a locally irreducible singular
point of $P$. If $s$ is an ordinary node of $P$
%(i.e. an $A_1$-point)
then $b_s(P)=1$.
{\em The nodal number $b(P)$ is the sum of the local nodal numbers
over all singular points of $P$.}

In many applications the numbers $c(P)$ and $b(P)$ play the r\^ole
of the total numbers of cusps and nodes even if $P$ has higher singularities.

\subsection{Integration w.r.t.\ Euler characteristic}\label{Viro-calculus}

Because of the equation
$$
\chi(A \cap B) = \chi(A) + \chi(B) - \chi(A\cap B)
$$
the Euler characteristic may be viewed as a signed
finite-additive measure. If $f$ is a function on
$X$ taking finite many values, its integral
with respect to the Euler characteristic is
defined by
$$
\int_X f\, d\chi  = \sum_{y}  y \, \chi(f^{-1}(y)) \,.
$$
Calculus of such integrals was developed by O.~Viro
\cite{V}. Applications of this calculus to the classical
Pl\"ucker and Klein formulas may also be found in \cite{V}.

More generally, given a map
$$
F: X \to Y
$$
one may push-forward $F_*f$ of  $f$ by the formula
$$
\left[F_*f\right](y) = \int_{F^{-1}(y)} f \, d\chi \,.
$$
Under additional hypotheses on $F$ and $f$, this operation has
natural functorial properties like $(FG)_*=F_* G_*$, see e.g.\
Section 7.3 in \cite{W}. The case of main importance for
us will be when $F$ is a nonconstant map of smooth curves,
in which case,
$$
F_*1 = \deg F - \br F \,,
$$
where $\br F$ is the \emph{branch divisor} of $F$, that is,
the sum of all critical points of $F$ with their
multiplicities.

\section{The equation of $R$}

\subsection{$R$ as a Resultant}\label{sRhat}

Consider the curve $\cR$ defined by the equations
\begin{equation}
  \label{eqRh}
     P(az,bw)=Q(z,w)=
\det
\begin{pmatrix}
\frac{\partial}{\partial z} P(az,bw) &  \frac{\partial}{\partial z} Q
\vspace{2 mm}\\
\frac{\partial}{\partial w} P(az,bw) &  \frac{\partial}{\partial w} Q
\end{pmatrix} = 0 \,.
\end{equation}
in $\tor\times\tor$. We have $\cR\subset\cS$, where
the surface $\cS$ was defined in \eqref{cS}. In general,
$\cR$ may have several
components. The Wronskian in \eqref{eqRh} vanishes at
any singular point of $q$ of $Q$, hence such a point contributes
a translate of  $P$ to $\cR$. The multiplicity of
this component equals $m(q)+\mu(q)-1$, see the proof of
Proposition \ref{dGauss}. Symmetrically, singular points of
$P$ contribute copies of $Q$ to $\cR$. The remaining
components correspond to actual tangency. The curve $R$
is the projection of these components to the $(a,b)$-plane.

The equation of $R$ may be found by eliminating
$z$ and $w$ from \eqref{eqRh}, followed by factoring off
spurious components caused by singularities. Gr\"obner basis
algorithms give one way to perform this elimination.
We find that in practice it is easier and faster to
use resultants for this computation.

Consider the polynomial
$$
R_1(a,b,z)=\textup{resultant}_w\big(P(az,bw),Q(z,w)\big) \,.
$$
The equation $R_1(a,b,z)=0$ cuts out the image $\cS_1$ of
$\cS$ under the projection along the $z$ direction.
This projection may create singularities, specifically
$\cS_1$ may intersect itself along curves. Consider
$$
R_2(a,b)=\textup{resultant}_z\left(R_1,
\frac{\partial}{\partial z} R_1\right) \,.
$$
Its zero locus consists of critical points of the
projection of $\cS_1$ along the $z$ direction, together
with contributions of singularities. Those can be
identified and removed by factoring $R_2$ as they
all occur with multiplicity greater than one. For
example, the image of double point curves of $\cS_1$
will occur in $R_2$ with multiplicity two.

Exact factorization is possible and rather effective
for polynomials with rational coefficients. In
probabilistic applications, the coefficients may
be known only approximately, making exact
factorization impossible. This is why it is useful
to know the Newton polygon of $R$, which will be
determined below. Also note that sometimes, albeit
rarely, components of $R$ may occur with
multiplicity in $R_2$ as, for example, in \eqref{double_dual}.

\subsection{Multiplicity of tangency}

Given two curves $C$ and $C'$ and a point $p$
define their tangency multiplicity at $p$ as
$$
(C \curlywedge C')_p = (C \cdot C')_p - m_{C}(p)\,  m_{C'}(p)
$$
where $(C \cdot C')_p$ is the local intersection
multiplicity of $C$ and $C'$ at $p$ and $m_C(p)$ is
the multiplicity of the point $p$ on $C$. Note that
$(C \curlywedge C')_p =0$ if $p\notin C$ or if $C$ and $C'$
intersect transversely at $p$. Also note that the
tangency multiplicity is additive over branches of both
$C$ and $C'$.

We have the following

\begin{prop} The multiplicity of a point $(a,b)\in R$ is
the total tangency multiplicity of $P(az,bw)=0$ and
$Q(z,w)=0$, assuming they are not tangent at
infinity.
\end{prop}

In other words, the local cuspidal number at each branch
of a singular point of $R$ is equal to the corresponding tangency
multiplicity minus one.

\begin{proof}
Since multiplicity is a local quantity, we can work with
the additive version \eqref{defadd}.
We may work with each branch of the singularity individually.
Suppose $P$ and $Q$
are tangent at $z=w=0$ and that $w=0$ is their
common tangent. In this case the origin will belong to
$P-Q$ and $b=0$ will be the corresponding tangent.
We need to compute the total multiplicity $m_R$ of all
branches of $R$ corresponding to this tangency. It
equals the number of
roots $b_i$ of $R(\epsilon,b)=0$ with $b_i\approx 0$
where $\epsilon$ is a fixed small nonzero number. For $b=b_i$
the curves
\begin{equation}
  \label{PQeps}
P(z+\epsilon,w+b)=0\,, \quad  Q(z,w)=0\,,
\end{equation}
are tangent.

Consider the intersection $U$ of the curve defined by
\eqref{PQeps} in the $(z,w,b)$-space with a neighborhood
of the origin. We will compute the Euler characteristic
of $U$ in two different ways. Viewing $U$ as a degree
$m_{P}\,  m_{Q}$ branched
covering of the $z$ line, we get
$$
\chi(U) = - m_{P}\,  m_{Q} + m_{P} +  m_{Q} \,,
$$
where $m_P$ and $m_Q$ is the multiplicity of
the origin on $P$ and $Q$, respectively.

On the other hand, we may view $U$ as a branched
covering of the $b$-line. It has degree $(P\cdot Q)_0$,
where $(P\cdot Q)_0$ is the intersection multiplicity of
$P$ and $Q$ at the origin. We have
$$
\chi(U) = (P\cdot Q)_0 - m_P (m_Q-1) - m_Q (m_P-1) - m_R
$$
because each tangency in \eqref{PQeps} corresponds to
a simple branchpoint. This concludes proof.
\end{proof}

\section{Newton polygon of log-front}

\subsection{Formula for $\Delta_R$}

Let $\Delta_P$ and $\Delta_Q$ denote the marked Newton
polytopes of $P$ and $Q$. We denote by $-\Delta_Q$ the
reflection of $\Delta_Q$ about the origin, with the
corresponding marking.

We say that two edges $E$ and $F$ are opposite
and
write $E\apar F$, if their outward normals point in the opposite
direction. As we will see, in the simplest case when
no edge of $-\Delta_Q$ is opposite to an edge of $\Delta_P$,
the Newton polygon $\Delta_R$
is simply the Minkowski sum
\begin{equation}
  \label{Minks}
    \deg(\gamma_Q) \cdot \Delta_P  + \deg(\gamma_P) \cdot (-\Delta_Q)\,.
\end{equation}
Every pair of opposite edges
$E\subset \Delta_P$ and $F\subset -\Delta_Q$
makes the polygon $\Delta_R$ smaller than \eqref{Minks}.

Given two partitions $\lambda$ and $\mu$, we define
\begin{align}
\langle \lambda,\mu \rangle  & = \sum_i \lambda'_i \, \mu'_i
\label{lan1} \\
& = \sum_{i,j} \min(\lambda_i,\mu_j) \label{lan2} \,,
\end{align}
where $\lambda'$ denotes the conjugate partition.
The equivalence of \eqref{lan1} and \eqref{lan2} is
an elementary combinatorial fact.

Given an edge $E$, let
$$
\mathbf{n}_E = \frac{E}{|E|}
$$
is the unit (non-oriented) lattice segment in the direction of $E$.
Note that (similarly to Newton polygons) such lattice segments
are defined only up to translations by $\Z^2$.

\begin{thm}\label{thmNewton}
The Newton polygon of the log-front $R$ is
the unique polygon $\Delta_R$ satisfying
\begin{equation}
  \label{polRes}
\Delta_R + \sum_{E\apar F}
\langle \lambda(E),\lambda(F) \rangle \, \mathbf{n}_E=
\deg(\gamma_Q) \cdot \Delta_P  + \deg(\gamma_P) \cdot (-\Delta_Q)\,,
\end{equation}
where the summation is over all pairs $E\subset \Delta_P$ and
$F\subset  -\Delta_Q$
of opposite edges and  $\lambda(E)$ are the corresponding partition
markings, i.e.\ multiplicities of intersection with the boundary.
\end{thm}

Minkowski ``subtraction'' of segments, implicit in formula \eqref{polRes},
simply reduces the length of all edges in that direction by the
given amount.

\subsection{Proof}\label{proofPolygon}

\begin{figure}[!htbp]
  \centering
  \scalebox{0.5}{\includegraphics{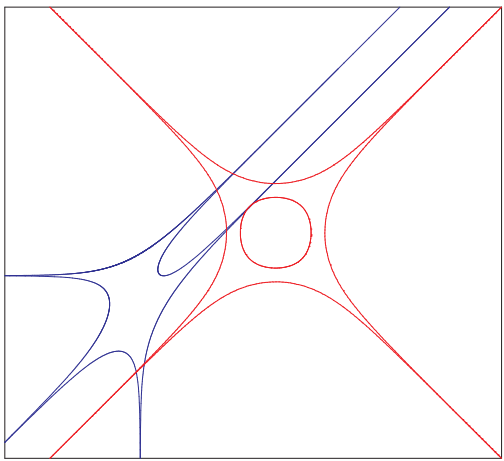}}
  \caption{Tangency to a tentacle creates a boundary
point of a log-front}
  \label{fig_tang_amoeb}
\end{figure}

By definition, a boundary point of $R$ is a tangency of
the curves $P$ and $Q$, dilated with respect to one
another by an infinite amount. In logarithmic coordinates,
infinite dilation becomes an infinite shift. Shifted
by an infinite amount, the amoeba $\Log(P)$ of $P$
is either a union of parallel lines (tentacles) or
empty, in which case there will not be any
boundary points of $R$ in this direction. This is
illustrated in Figure \ref{fig_tang_amoeb}.

Tangency to a tentacle occurs at points mapped by
the logarithmic Gau\ss\ map to the slope of the tentacle.
Therefore, if $\Log(P)$ and $\Log(Q)$ do not have
tentacles in the same direction (that is,
$\Delta_P$ and $\Delta_Q$ have no parallel edges) then
each tentacle of $\Log(P)$ contributes $\deg \gamma_Q$
points to $\partial R$. The multiplicity of the resulting point
in $\partial R$ equals
the multiplicity of the corresponding point of $\partial P$.

The effect of a tentacle-to-tentacle tangency may be
studied in a local model, for example,
$$
P=\{w-(z-1)^k\}\,, \quad Q=\{w-(z-1)^l\}\,.
$$
In this case,
$$
\gamma_P(z)\sim\frac{k}{z-1}\,, \quad
\gamma_Q(z)\sim\frac{l}{z-1}\,, \quad z\to 1 \,,
$$
Equating $\gamma_P(az)=\gamma_Q(z)$ gives a parametrization
of $R$ such that
\begin{equation}
  \label{ttt}
  a \sim 1, \quad b\propto (z-1)^{k-l}\,, \quad
 z\to 1 \,.
\end{equation}
This means that tentacles of $P$ and $Q$ pointing in the
opposite direction (which means that $k$ and $l$ have opposite
sign) produce a tentacle of combined multiplicity
$|k|+|l|$ pointing in the same direction as
the $P$-tentacle. Tentacles pointing in the same direction,
by contrast,  produce a tentacle of multiplicity
$|k-l|$ pointing in one of the two directions.
This rule can be phrased more naturally in terms of the
reflected Newton polygon $-\Delta_Q$ because the
reflection flips the direction of tentacles.

The polygons $\Delta_P$ and $-\Delta_Q$ can, in total,
have as many as 4 edges with the same slope.
Let $E,\barE\subset\Delta_P$
and $F,\barF\subset-\Delta_Q$ be such a 4-tuple of
edges and let $\lambda(E),\dots$ denote the
partition marking.
Assume that the outward normals to $E,F$ point
in the same direction. By our computations, the polygon
$\Delta_R$ will have an edge in the same direction as $E$ and $F$
of length
\begin{gather}
|E| (\deg \gamma_Q - \ell(\lambda(F)) - \ell(\lambda(\barF))) +
|F| (\deg \gamma_P - \ell(\lambda(E)) - \ell(\lambda(\barE)))
\notag \\
+\sum_{i,j} \left(\lambda(E)_i +\lambda(F)_j\right) +   \label{|G|}\\
\sum_{i,j} \max(\lambda(E)_i-\lambda(\barF)_j,0) +
\sum_{i,j} \max(\lambda(F)_i-\lambda(\barE)_j,0) \notag \,.
\end{gather}
Note that the middle line in \eqref{|G|} cancels with a part of
the first line. A further cancellation is obtained from
$$
- |E| \, \ell(\lambda(\barF)) +
\sum_{i,j} \max(\lambda(E)_i-\lambda(\barF)_j,0) =
- \sum_{i,j} \min(\lambda(E)_i,\lambda(\barF)_j)
$$
using formula \eqref{lan2}. This concludes proof.

Counting the boundary points of $R$ without respect to the
corresponding direction $\mathbf{n}_E$ gives us the following corollary.

\begin{coro}
The cardinality of $|\dd R|$ equals to
$$\deg\gamma_P|\dd\Delta_Q|+\gamma_Q|\dd\Delta_P|-\sum_{E\apar F}
\langle \lambda(E),\lambda(F) \rangle.$$
%Recall that $|\dd\Delta_P|$ is the number of lattice points on the boundary
%of $\Delta_P$.
\end{coro}

\subsection{Examples}\label{expl_poly}

\begin{figure}[!htbp]
  \centering
  \scalebox{0.4}{\includegraphics{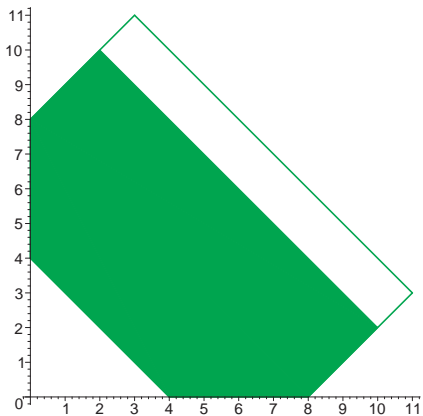}}
  \caption{Possible $\Delta_R$ for $\Delta_P$ and $\Delta_Q$
from Figure \ref{fig_Polygons}.}
  \label{fig_newton_front}
\end{figure}

Let $P$ and $Q$ be generic with Newton polygons
from Figure \ref{fig_Polygons}. In this case,
the Newton polygon $\Delta_R$ is the larger of the
two polygons plotted in Figure \ref{fig_newton_front}.
If $P$ develops a tangency to the boundary then
$\Delta_R$ shrinks and becomes the smaller polygon
in Figure \ref{fig_newton_front}. What happens
in this case, is that $R$ becomes reducible with
one component being a boundary divisor. The boundary divisor
corresponds to tangency with $P$ at the newly
developed point of tangency to infinity.

As another example, consider the case when $P$ and $Q$
are generic curves of degree $d_P$ and $d_Q$.
In this case, both $\Delta_P$ and $\Delta_Q$ are triangles
with sides of length $d_P$ and $d_Q$, while
$\Delta_R$ is a hexagon with sides of length
$$
  d_P d_Q (d_Q-1), d_Q d_P (d_P -1), \dots\,,
$$
cyclically repeated. In particular, when $d_P=1$ we get
a triangle with side $d_Q (d_Q-1)$, reproducing
the very classical formula
for the degree of dual curve.

In the example in Figure \ref{froz}, the degrees of
logarithmic Gauss\ maps of $P$ and $Q$ are $4$ and $10$
respectively, hence $\Delta_R$ is a hexagon
with sides $8$ and $12$.

\subsection{Reconstructing the curve $Q$ by the log-front $R$}
It is instructive to have a closer look at the
computation done in \eqref{ttt} in the case
$k=l$. In this case, the point $(a,b)$ does
not escape at infinity, but it is still a pole
of the logarithmic Gauss map of $R$ (see also
Section \ref{sGR} below). In other words,
the curve $R$ is tangent to
$a=1$ at the corresponding point.

In probabilistic applications,
the curve $P$ is given and one knows that the frozen
boundary is compact and tangent to given
lines. This allows to fix the real boundary
points of $Q$ and their multiplicities
(which have to match the corresponding
multiplicities for $P$). In particular,
this is how the curve $Q$ in Figure
\ref{PandQ} is determined from the
requirement that the frozen boundary
in Figure \ref{froz} is inscribed in
a hexagon.

\section{Pl\"ucker-type formulas for $R$}

\subsection{Logarithmic Gau\ss\ map of $R$}\label{sGR}

Let $\Rh$ be the product
of the normalizations $\Pt$ and $\Qt$ over their
logarithmic Gau\ss\ maps:
\begin{equation} \label{Gauss_diag}
\xymatrix{
& \Rh \ar[dr]^{\pi_Q} \ar[dl]_{\pi_P}  \ar[dd]^{\gamma_R}\\
\Pt \ar[dr]_{\gamma_P}&& \Qt \ar[dl]^{\gamma_Q}
  \\
& \pline
}
\end{equation}
In English, a point of $\Rh$ corresponds to a tangency of
a branch of $P$ to a translate of a branch of $Q$.
Note that $\Rh$ is a partial normalization of the plane curve $R$:
we have the factorization $\Rt\to\Rh\to R$ for the normalization $\Rt\to R$.
We claim that the logarithmic Gau\ss\ map of
$\Rt$ factors through the natural map $\Rt\to \Rh$.

\begin{prop}
The logarithmic Gau\ss\ map for $R$ is the composition
$$
\Rt \to \Rh \overset{\gamma_R}\to \pline
$$
\end{prop}

\begin{proof}
The essential geometric content of this result is
already implicit in \eqref{double_dual}. In
coordinates, the claim is elementary to check on the open dense set
where tangency is
nondegenerate.  It suffices to analyze the
Gauss\ map of the additive analog $P-Q$ of the
log front $P\div Q$. Let $P$ and $Q$ be given by
$$
y=P(x)\,,\quad y=Q(x) \,.
$$
Then $P-Q$ is parametrized by
$$
(x-s,P(x)-Q(s))
$$
where $s=s(x)$ is a solution of $P'(x)=Q'(s)$. The
Gauss\ map of the above curve is clearly
$P'(x)=Q'(s)$.
\end{proof}

\begin{cor}\label{T_deg}
\begin{equation}
  \label{gamma_R}
\deg \gamma_R = \deg \gamma_P \, \deg \gamma_Q \,.
\end{equation}
\end{cor}

\subsection{Geometric genus of $R$}

The geometric genus of $R$, or equivalently, the
Euler characteristic $\chi(\Rt)$ can be computed by
Riemann-Hurwitz formula applied to the
map $\gamma_R$ in \eqref{Gauss_diag}. To do this we
need to be able to compare $\chi(\Rt)$ and $\chi(\Rh)$.

The map $\gamma_R$ is ramified over the branchpoints
of either $\gamma_P$ or $\gamma_Q$. Suppose that
at a point $p\in\Pt$ the map $\gamma_P$ is given by
\begin{equation}
  \label{coordx}
    x \to x^{\nu}\,, \quad \nu=1,2,\dots\,,
\end{equation}
in a suitable local coordinate $x\in \C$ centered at $p$.
In other words, suppose that $p$ is a logarithmic
inflection point of multiplicity $\nu-1$.
Note that such point may be singular or non-singular point of $P$.
This
integer $\nu$ will be denoted by $\nu(p)$.

A point of $\Rh$ is a pair $(p,q)\in \Pt \times \Qt$ such that
$$
\gamma_P(p) = \gamma_Q(q) = x \in \pline \,.
$$
Let $x_p$ and $x_q$ be local coordinates at $p$ and $q$
as in \eqref{coordx}. The local equation of $\Rh$
is
$$
x_p^{\nu(p)} = x_q^{\nu(q)} \,, \quad (x_p,x_q)\in \C^2 \,.
$$
This has $\gcd(\nu(p),\nu(q))$ branches and hence
produces $\gcd(\nu(p),\nu(q))$ points
$$
r_1,\dots, r_{\gcd(\nu(p),\nu(q))} \in \Rt
$$
that are mapped to $x$ by the logarithmic Gau\ss\ map.
 Each is a logarithmic inflection point
of multiplicity $\lcm(\nu(p),\nu(q))-1$. Note that while
all $r_i\in \Rt$ are logarithmic
inflection points of the same multiplicity, the
multiplicities of the corresponding branches of $R$ may
be different and are not determined by the numbers $\nu(p)$
and $\nu(q)$.

{}From definitions
\begin{equation}
  \label{chiRh}
    \chi(\Rh) = \int_{\Rh} 1 \, d\chi  = \int_{\pline} (\gamma_P)_* 1 \,
(\gamma_Q)_* 1 \, d\chi \,.
\end{equation}
We conclude

\begin{thm}\label{T_gen}
The Euler characteristic $\chi(\Rt)$ is given by
$$
\chi(\Rt) = \chi(\Rh) + \sum_{\gamma_P(p)=\gamma_Q(q)} (\gcd (\nu(p),\nu(q)) -1)
$$
where $ \chi(\Rh) = ((\gamma_P)_* 1 ,
(\gamma_Q)_* 1)_\chi =
- 2 \deg \gamma_P \, \deg \gamma_Q +
\chi(\Pt)\, \deg \gamma_Q + \chi(\Qt) \, \deg \gamma_P$.
\end{thm}
The above inner product with respect to the Euler characteristic
is defined as in \eqref{chiRh}. Note that generically the
the branchpoints of $\gamma_P$ are disjoint from branchpoints of $\gamma_Q$.
In this case,
the above formulas may be simplified as follows.
\begin{cor}
If tangency does not occur at two logarithmic inflection points then
\begin{equation}
  \label{chiRs}
  \chi(\Rt) = \chi(\Rh) = - 2 \deg \gamma_P \, \deg \gamma_Q +
\chi(\Pt)\, \deg \gamma_Q + \chi(\Qt) \, \deg \gamma_P  \,.
\end{equation}
\end{cor}

\subsection{Nodes and cusps of $R$}

%\ignore{
The total cuspidal number of $R$ may now be
determined using the formula \eqref{n_cusps}.

\begin{thm}\label{cuspidal}
We have the following expressions for the cuspidal number $c_{{R}}$
of the resultant curve $R=R(P,Q)$.
$$c(R)=-\chi(\Rt)+|\dd R|
-(\chi(\Pt)+c(P)-|\dd P|)(\chi(\Qt)+c(Q)-|\dd Q|).$$
\end{thm}

\begin{proof}%[Proof of Theorem \ref{cuspidal}]
This formula can be obtained as a straightforward
combination of \eqref{n_cusps} and Corollary \ref{T_deg}.
%and Theorem \ref{T_gen}.
Nevertheless it is instructive to prove it by the Viro
calculus (see Section \ref{Viro-calculus}) to prepare
a way for the real counterpart in the next section.
%(note that \cite{V} already contains application of this
%calculus to the Pl\"ucker and Klein formulas in the classical case).

%Recall that $\Delta_{PQ}=\Delta_Q-\Delta_P$ was defined as the Minkowski
%difference.
Consider the family of translates $\tau_s({P})$ parameterized
by $s\in\tor$.
%$s\in\cT\setminus\ver(\cT)$,
%where $\ver(\cT)$ are the vertices of the toric surface
%$\cT$ corresponding to the polygon $\Delta_{PQ}$ as in \eqref{cT}.
%i.e. the fixed points of the $\tor$-action.
%Clearly, $\chi(\cT\setminus\ver(\cT))=0$.
%
Denote with $X(s)$
%$\#(\tau_s(\bar{P}),\bar{Q})$
the Euler characteristic of
the space of pairs $(t_P,t_Q)\in\tilde{P}\times\tilde{Q}$,
such that $\tau_s(t_P)=t_Q$.
%and either $t_P\notin\dd P$ or $t_Q\notin\dd Q$.
%$\tau_s(\nu_{\bar{P}}(t_P))=\nu_{\bar{Q}}(t_Q)$
%and either $\nu_{\bar{P}}(t_P)\in\tor$ or $\nu_{\bar{Q}}(t_Q)\in\tor$.
We have
\begin{equation}\label{cint}
\int\limits_{\tor} X(s) d\chi(s)=
(\chi(\tilde{P})-|\dd P|)(\chi(\tilde{Q})-|\dd Q|).
%(l_{\dd\Delta_P}-\sigma_P)(l_{\dd\Delta_Q}-\sigma_Q).
\end{equation}
It may be viewed as a corollary of the Fubini theorem
since each pair of points from $\tilde{P}\setminus\dd P$ and
$\tilde{Q}\setminus\dd Q$ will
appear in $X(s)$ once.
%, unless both of them
%are mapped to the boundary of $\cT$.
On the other hand we have
\begin{multline}
\label{cT2}
\int\limits_{\tor} X(s) d\chi(s)=
%l_{\dd\Delta_P}l_{\dd\Delta_P}-\delta(P,Q)
-\chi(\tilde{R})+|\dd R|\\-c(P)(\chi(\tilde{Q})-|\dd Q|)
-c(Q)(\chi(\tilde{P})-|\dd P|)
-c(P)c(Q)-c(R).
\end{multline}
Indeed, we have $\chi(\tor)=0$ while
$X(s)$ is constant for generic $s\in\cT$
(namely, by Bernstein-Kouchnirenko formula it is equal to
$\Vol(\Delta_P,\Delta_Q)$).
However the value $X(s)$ drops
if $s\in R$ and it drops further (by $m-1$) if there is a singular point
of multiplicity $m$ at a branch of $\tilde{R}\setminus\dd R$ through $s$
or if $\tau_s(\Pt)$ and $\Qt$ have a tangency of higher order
(which in turn corresponds to a cusp of $\tilde{R}$).
%As $\int\limits_{\C T_{\Delta}\setminus\ver(\C T_{\Delta})}
%\#(\tau_s(\bar{P}),\bar{Q}) d\chi(s)$
%also compute some pairs $(t_P,t_Q)$ with
%$\nu_{\bar{P}}(t_P),\nu_{\bar{Q}}(t_Q)\notin\tor$
%we have to insert the correction by $l_{\dd\Delta_P}l_{\dd\Delta_P}-\delta(P,Q)$.
We get the theorem as the combination of \eqref{cint} and \eqref{cT2}.
\end{proof}

If $P$ and $Q$ are generic, %(in a proper sense),
the only singularities of $R$ will be cusps and nodes and the
nodal number of $R$ may be recovered from the adjunction formula.

For example, let $P$ and $Q$ be generic curves
of degree $d_P$ and $d_Q$. In this case
%we have $\bar\delta(P,Q)=3d_Pd_Q$,
%%$l_{\dd\Delta_P}=3d_P$, $l_{\dd\Delta_Q}=3d_Q$,
%$\chi(\tilde{P})=3d_P-d_P^2$, $\chi(\tilde{Q})=3d_P-d_Q^2$
%and
$$
\chi(\Rt) = -d_P d_Q (4 d_P d_Q -3 d_P-3 d_Q)\,.
$$
Note that this grows very fast with $d_P$ and $d_Q$.
Already when $P$ and $Q$ are generic conics, the
geometric genus of $R$ is $9$.

The Newton polygon
of $R$ and its boundary  were determined above in Section \ref{expl_poly}.
Generically, $R$ will not  have multiple point
on the boundary. For the number of cusps Theorem \ref{cuspidal} produces
$$
\textup{\# cusps} = 3 d_P^2 d_Q^2-6 d_P d_Q \,.
$$
Accordingly for the number of nodes we
get:
$$
\textup{\# nodes} =
\frac{d_P^4 d_Q^2}{2}+2 d_P^3 d_Q^3+\frac
{d_P^2 d_Q^4}{2}-3 d_P^3 d_Q^2-3 d_P^2 d_Q^3
-2 d_P^2 d_Q^2+9 d_P d_Q\,.
$$
For $d_P=1$, these specialize to the classical
Pl\"ucker formulas.
The last two formulas were
obtained in \cite{FO}.

\section{Klein-type formula for the log-front $R$}
\subsection{Refinement of the nodal and cuspidal numbers}
For this section it will be important that both $P$ and $Q$
are defined over the field $\R$ of real numbers. In this case, clearly,
the curve $R$ is defined over $\R$ as well.

An algebraic curve $P\subset\cT$ is defined over $\R$
if and only if it is invariant with respect to the involution
of complex conjugation $\conj:\cT\to\cT$.
The fixed point locus of this involution is the real
toric surface $\R\cT\subset\cT$.
The real locus $\R{P}$ coincides with the intersection
$P\cap\R\cT$.

For a real curve $P$ we may refine both the cuspidal number $c(P)$
and the nodal number $b(P)$ as follows.
Let $c^{\re}(P)$ be the number of cusps of $P$
(counted with multiplicity as in Section \ref{bcdefn}) that are
%from the real branches of singularities of $P$.
{\em real}, i.e. contained in $\R T_\Delta$.
In other words, to get $c^{\re}(P)$ we add over (real and imaginary)
branches of {\em real} singular points of $P$ the multiplicities of these
branches diminished by 1.
We set $c^{\im}(P)=c(P)-c^{\re}(P).$

Let $b^{\re}_+(P)$ be the sum over all singular
points $s$ of $P$ in $\R\cT$ of the number
of the pairs of conjugate imaginary branches of $P$.
%Equivalently, $b^{\re}_+(P)$ is one half of the
%cardinality of $\nu_{\bar{A}}^{-1}(\R T_{\Delta_A})$
%(recall that $\nu_{\bar{A}}:\tilde{A}\to \bar{A}$ is the
%normalization).
If $s\in\R P$ is a singular point
we may introduce $b^{\re}_-(P,s)$ as the sum
of the local intersection numbers over all possible
pairs of the real branches
of $P$ at $s$. Let $b^{\re}_-(P)$
be the sum of $b^{\re}_-(P,s)$ over all
real singular points of $P$ and let
$b^{\im}(P)=b(P)-b^{\re}_+(P)-b^{\re}_-(P)$.

Thus we get the refinements $$c(P)=c^{\re}(P)+c^{\im}(P)$$
and $$b(P)=b^{\re}_+(P)+b^{\re}_-(P)+b^{\im}(P).$$
In the case when the only singularities of $P$ are ordinary cusps and
nodes the numbers $c^{\re}(P)$ and $c^{\im}(P)$
are the numbers of real and complex cusps respectively while
the numbers $b^{\re}_-(P),b^{\re}_+(P)$ and $b^{\im}(P)$
are the numbers of real hyperbolic ($A_1^-$), real elliptic ($A_1^+$)
and imaginary nodes.

It is convenient to define the {\em boundary nodal number}
$b(\dd P)$ to be equal to the sum of the multiplicities
of the boundary points of $P$ minus the simple cardinality of $\dd P$.
(Note that since $\dd P$ is 0-dimensional this boundary nodal number
also works as the boundary counterpart of the cuspidal number.)
Again we have the refinement $b(\dd P)=b^{\re}(\dd P)+b^{\im}(\dd P)$
that counts real and imaginary boundary points separately.
Thus $b^{\re}(\dd P)$ is the measure of nontransversality of the real
locus $\R P$ to the boundary divisor of the toric surface $\cT$.

\subsection{Computations for $R$}

In general, the set of refined nodal and cuspidal numbers
$b^{\re}_+,b^{\re}_-,c^{\re},c^{\im}$ for
the log-front $R$ is not determined by
the corresponding sets for $P$ and ${Q}$.
Let us recall a well-known example illustrating this
in the case of classical projective duality, i.e. if
$P$ is given by the polynomial $z+w+1$: two real inflection
points may disappear together with a bitangent real line.
\begin{exa}
Let $\R Q_1,\R Q_2\subset\rp^2$ be real quartic curves
pictured in Figure \ref{quart}. The first curve can be constructed
by perturbation of the union of two ellipses while the second one
can be constructed by perturbation of the union of four lines.
We have $b^{\re}_+({R_1})=0,c^{\re}({R_1})=8,
b^{\re}_+({R_2})=4,c^{\re}({R_1})=0$.
\begin{figure}[h]
\centerline{\psfig{figure=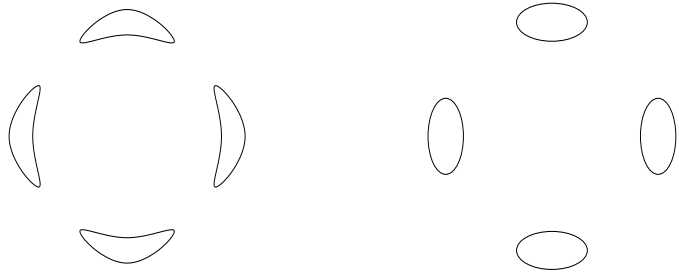,height=1.4in,width=3.5in}}
\caption{\label{quart} Real quartic curves $\R Q_1$ and $\R Q_2$
whose projective duals have distinct $b^{\re}_+$ and $c^{\re}$.}
\end{figure}
\end{exa}

The classical Klein's formula \cite{Kl-paper}
allows us to compute
$2b^{\re}_+(\bar{R})+c^{\re}(\bar{R})$
in the case when $P$ is a line.
To generalize this statement for a larger class of real curves $P$
let us look at the
%so-called {\em coamoebas} of the curves
%$P,Q\subset\tor$.
argument map $\Arg:\tor\to S^1\times S^1$
defined by $$\Arg(z,w)=(\arg(z),\arg(w)).$$
The image $\Arg(A)\subset S^1\times S^1$
is called {\em the coamoeba} or {\em the alga} of $A\subset\tor$, cf. \cite{Pa},
\cite{KeVa}.

The real 2-torus $S^1\times S^1$ is a group which has $\Z_2\times\Z_2$
as its subgroup.
Let $T=S^1\times S^1/\Z_2\times\Z_2$ be the quotient group and
$\beta:S^1\times S^1\to T=S^1\times S^1/\Z_2\times\Z_2$ be
the projection map. Note that $T$ is itself a group isomorphic
to $S^1\times S^1$ and the zero in this group is $0_T=\beta(1,1)$.

We need to compactify the map $\beta\circ\Arg:P\setminus\dd P\to T$.
Note that this map does not extend to $P$.
Let us consider $\hat{P}$ to be {\em the real blow-up} of $\tilde{P}$
at the finite collection of points
$\dd P\subset\tilde{P}$.
Naturally, $\hat{P}$ is a closed non-orientable surface whose
Euler characteristic coincides with that of $\tilde{P}\setminus{\dd P}$
(in the case of non-compact spaces
we use Euler characteristic for homology with closed support)
Furthermore, we have a natural extension
$$\Alga_P:\hat{P}\to T$$
such that $\Alga_P|_{\tilde{P}\setminus\dd P}=
\beta\circ\Arg|_{\tilde{P}\setminus\dd P}$.

%Note that the involution of complex conjugation $\conj:\bar{A}\to\bar{A}$
%of a curve $\bar{A}$ lifts to its normalization $\tilde{A}$.
%We may define $\R\tilde{A}$ as the fixed point set of this lift
%so that $\nu_{\bar{A}}(\R\tilde{A})\subset\R\bar{A}$.
%Define the function
%$$\rho:\tilde{Q}\to\Z$$
%to be $\rho(s)=\chi((\hat\Arg_P)^{-1}(\Arg_Q(s)))$,
%if $s\in Q$ and
%$\rho(s)=...$

Consider the following integral
$$I_{\R}=\int\limits_{t\in\R\cT\setminus\ver(\R\cT)}
X(t) d\chi_{\R\cT}(t)$$
(note that here we use the family $\tau_t(P)$
only for $t\in\R\cT$).
This integral is the real counterpart of the integral \eqref{cint}
from the proof of Theorem \ref{cuspidal} where we defined the number
$X(t)$, but now we use only the real translations.
%was defined, but is defined
%using only real translations.

As in Theorem \ref{cuspidal} the integral $I_{\R}$ can be computed
in two different ways.

For this computation to depend only on visible characteristics
of the curves $P$ and $Q$ we need to assume that $P$ is a
Harnack curve (see Appendix A).
%The right-hand side of
%%Proposition \ref{generalKlein} contains the quantity
%\eqref{1I}
%may not be immediately visible from the geometry of $\R P$ and $\R Q$.
%However, if $P$ is a simple Harnack curve (see eg. \cite{MR})
%then we can compute this quantity explicitly with the help of
%the formulas deduced in Appendix A.

Let $$\epsilon(\Delta_P,\Delta_Q)=\sum\limits_{\Delta'_P,\Delta'_Q}
\Area(\Delta'_P+\Delta'_Q),$$
where $\Delta'_P$ (resp. $\Delta'_Q$) run over all possible
sides of the polygon $\Delta_P$ (resp. $\Delta_Q$) and
$\Delta'_P+\Delta'_Q$ stands for the parallelogram obtained
as the Minkowski sum of the intervals $\Delta'_P$ and $\Delta'_Q$.

%Let $\perimeter(\Delta_Q)$ be the number of lattice point on
%the perimeter of $\Delta_Q$, i.e. $\#(\Delta_Q\cap\Z^2)$.
\begin{prop}\label{intHarnack}
If $P$ is a Harnack curve and $Q$ is any curve defined over $\R$
then
%\begin{multline*}
%\int\limits_{\alpha\in T}
%\chi((\hat\Arg_P)^{-1}(\alpha))\chi((\hat\Arg_Q)^{-1}(\alpha))d\chi_T(\alpha)\\
%-\int\limits_{\alpha\in T}
%\chi((\hat\Arg_P)^{-1}(\alpha)\cap\dd\tilde{P})
%\chi((\hat\Arg_Q)^{-1}(\alpha)\cap\dd\tilde{Q})
%\chi_T(\alpha)\\
\begin{multline*}
I_{\R}=2\Area(\Delta_P)(\chi(\tilde{Q})-2b^{\re}_+({Q})-b^{\re}(\dd Q)
-|\dd\Delta_Q|)+\\
(2b^{\re}_+({P})+b^{\re}(\dd P))(2b^{\re}_+({Q})+b^{\re}(\dd Q))-\epsilon(\Delta_P,\Delta_Q),
\end{multline*}
%\end{multline*}
%where $\epsilon(\Delta_P,\Delta_Q)$ is defined to be
\end{prop}
\begin{proof}
Note that if two points $s_P\in\hat{P}$ and $s_Q\in\hat{Q}$ are different
by a translation by $t\in\R\cT$ then
$\Alga_P(s_P)=\Alga_Q(s_Q)$.
Conversely, if $\Alga_P(s_P)=\Alga_Q(s_Q)$ and
either $s_P$ or $s_Q$
is not in the exceptional divisor of the blowup map
%and the normalization
%(of $\hat{P}$ or $\hat{Q}$ respectively)
then there
exists $t\in\R\cT\setminus \ver(\R\cT)$
such that $\tau_t({P})$ and ${Q}$ have
$s_Q$ as one of their
intersection point. The Euler characteristic of
the space formed by the pairs $(s_P,s_Q)$ with
$\Alga_P(s_P)=\Alga_Q(s_Q)$
and such that $s_P$ and $s_Q$ are from the exceptional
divisors of the blowups $\hat{P}\to\tilde{P}$ and $\hat{Q}\to\tilde{Q}$ respectively
is $\epsilon(\Delta_P,\Delta_Q)$.

We have $\chi((\Alga_P)^{-1}(\alpha))=2\Area(\Delta_P)$
for any $\alpha\in T\setminus\{0_T\}$ by Lemma A\ref{arg-harnack}.
Note that
%$\chi((\hat\Arg_P)^{-1}(0_T))=2b^{\re}_+(Q))$.
$\chi((\Alga_Q)^{-1}(T\setminus \{0_T\}))=\chi(\hat{Q})-2b^{\re}_+({Q})-b^{\re}(\dd Q)=
\chi(\tilde{Q})-|\dd\Delta_Q|-2b^{\re}_+({Q})-b^{\re}(\dd Q)$.
The contribution of $\alpha=0_T$ to the left-hand side of
Proposition \ref{intHarnack} is
$(2b^{\re}_+({P})+b^{\re}(\dd P))(2b^{\re}_+({Q})+b^{\re}(\dd Q))$ since $\R\hat{P}$
is a closed 1-manifold and thus $\chi(\R\hat{P})=0$.
\end{proof}

Let $n$ be the number of vertices of the polygon
$\Delta_{PQ}=\Delta_Q-\Delta_P$.
\begin{prop}\label{2prop}
If $P$ is a Harnack curve and $Q$ is any curve defined over $\R$
immersed near the boundary of $\R\cT$ then
\begin{multline*}
I_{\R}=\Vol(\Delta_P,\Delta_Q)(4-2n)\\-2\Area(\Delta_P)c^{\im}({Q})-(2b^{\re}_+(P)+
b^{\re}(\dd P))c^{\re}(Q)-
(2b^{\re}_+({R})+c_{\re}({R})).
\end{multline*}
\end{prop}
\begin{proof}
The Euler characteristic of $\R\cT\setminus\ver(\R\cT)$
is $4-2n$. For a generic $t$ we have $X(t)=
\Vol(\Delta_P,\Delta_Q)$ by the Bernstein-Kouchnirenko formula \cite{Be}, \cite{Ku}.
This number gets decreased if $\tau_t({P})$ and ${Q}$ are tangent
or if one of their intersection point is singular for a branch of $\tilde{Q}$
(note that ${P}$ is an immersed smooth curve since it is Harnack).
The latter case contributes $-2\Area(\Delta_P)c^{\im}({Q})
-(2b^{\re}_+(P)+b^{\re}(\dd P))c^{\re}(Q)$.

If $\tau_t({P})$ and ${Q}$ are tangent at a point in
a non-real point then we have a bitangency since both ${P}$
and ${Q}$ are invariant with respect to the involution of
complex conjugation. This contributes $-2b^{\re}_+({R})$.
\ignore{
To compute the Euler characteristic of
the tangencies at real points
we note that for every point $s_q\in Q$ we have exactly $a_P$
points such that there exists $t\in\R T_{\Delta_R}\setminus\ver(\R T_{\Delta_R})$
with $\tau_t(\bar{P})$ tangent to $\bar{Q}$ at $s_Q$ since there are
$a_P$ points in $\R\bar{P}$ with the same image of the Logarithmic Gauss
map as $s_Q$ since $P$ is simple Harnack (see \cite{Mi}).
For a point $s_q\in S^{\dd}_Q$ we have $a_P$ minus the number of points
at $S^{\dd}_Q$ with the same image of the Gauss map as $s_Q$ plus the number
of those points from this set that $\tilde{P}$ has the same local intersection
number with the boundary divisor as $\tilde{Q}$ has at $s_Q$. This contributes
$\sigma(\R P,\R Q)$.
}
The tangencies at real points contribute $\chi(\R\tilde{R})=0$ plus
$c_{\re}({R})$ where the tangencies of $\tau_t({P})$ and ${Q}$
are of higher order.
\end{proof}

\subsection{Klein's formula}
Combining Proposition \ref{intHarnack} and \ref{2prop}
we get the following theorem.
\begin{thm}\label{Harthm}
If $P$ is a simple Harnack curve and $Q$ is any curve defined over $\R$
immersed near the boundary of $\R\cT$
then
\begin{multline*}
2b^{\re}_+({R})+c^{\re}({R})=
(4-2n)\Vol(\Delta_P,\Delta_Q)+\\
2\Area(\Delta_P)(|\dd\Delta_Q|-\chi(\tilde{Q})
+2b^{\re}_+({Q})+b^{\re}(\dd Q)-c^{\im}({Q}))\\
-(2b^{\re}_+({P})+b^{\re}(\dd P))(2b^{\re}_+({Q})+b^{\re}(\dd Q)+c^{\re}(Q))
+\epsilon(\Delta_P,\Delta_Q).
\end{multline*}
\end{thm}
%\begin{proof}
%Since $\bar{P}$ is smooth we have the set $\Sigma\tilde{P}$ empty.
%Since $\chi(\hat\Arg_P)^{-1}(\alpha)=\Area(\Delta_P)$
%for $\alpha\in T\setminus 0_T$
%we have
%$$\sum\limits_{s\in\Sigma\tilde{Q}\setminus\nu_{\bar{Q}}^{-1}(\R T_\Delta)}
%c(s)\chi(\hat\Arg_P)^{-1}(\hat\Arg_Q(s))=\Area(\Delta_P)c^{\im}(\bar{Q}).$$
%\end{proof}

\begin{coro}\label{H-coro}
Let $\R{P},\R{Q}\in\rp^2$ be two curves of degree $d_P$ and $d_Q$
(respectively) not passing via $(0:0:1),(0:1:0),(1:0:0)$.
Suppose that $P$ is a simple Harnack curve and ${Q}$
is a smooth curve. Then
$$
2b^{\re}_+({R})+c^{\re}({R})=
d^2_Pd^2_Q-2d_Pd_Q.
$$
\end{coro}
\begin{proof}
We have $n=6$, $\Area(\Delta_P)=\frac{d_P^2}{2}$, $\chi(\tilde{Q})=3d_Q-d^2_Q$,
$|\dd\Delta_Q|=3d_Q$ and $\epsilon(\Delta_P,\Delta_Q)=6d_Pd_Q$.
\end{proof}

Let us deduce the Klein formula \cite{Kl}
in its classical form from Theorem \ref{Harthm}
in the case when $P$ is a line $z+w+1=0$
and $\R Q$ is a curve of degree $d$ in $\rp^2$ not passing
via $(0:0:1),(0:1:0),(1:0:0)$. In this case $\Delta_Q$
is a triangle with vertices $(0,0)$, $(d,0)$ and $(0,d)$.
The line $P$ is a simple Harnack curve of degree 1 with $a_P=1$
and thus
we may apply Theorem \ref{Harthm}. The polygon $\Delta_R$
is a hexagon in this case, so $n=6$, while $\Vol(\Delta_P,\Delta_Q)=d$
and $\epsilon(\Delta_P,\Delta_Q)=6d$.

Note that the classical Klein formula computes
$2b^{\re}_+({R})+c^{\re}({R})$
for the dual curve in $\R P^2$ while our formula does it for
${R}\subset\R\cT$. The toric surface
$\R\cT$ is the result of the blowup of $\R P^2$ at
three points and, in general, such blowup might change the
characteristics $b^{\re}_+$ and $c^{\re}$.
Let us assume that ${Q}$ intersects the boundary
divisor $\cp^2\setminus\tor$ transversely, so that this
blowup is disjoint from ${Q}$.
Then $\chi(\hat{Q})=\chi(\tilde{Q})-3d$.

We get
\begin{multline*}
2b^{\re}_+({R})+c^{\re}({R})=-8d
+2b^{\re}_+({Q})-\chi(\tilde{Q})+3d-c^{\im}({Q})+6d\\
=d-c^{\im}({Q})+2b^{\re}_+({Q})-\chi(\tilde{Q}).
\end{multline*}
Note that $c^{\im}({Q})=c_{{Q}}-c^{\re}({Q})$
while the classical Pl\"ucker formula \cite{Pl}
for ${Q}$ implies that
$\chi(\tilde{Q})+c_{{Q}}=2d-d^*,$
where $d^*$ is the degree of $R\subset\tor\subset\cp^2$
(called {\em the class of $Q$} in Klein's paper \cite{Kl-paper}).
Thus we get
$2b^{\re}_+({R})+c^{\re}({R})=
-d+d^*+2b^{\re}_+({Q})+c^{\re}({Q})$ and, therefore,
\begin{equation}
d-2b^{\re}_+({Q})-c^{\re}({Q})=
d^*-2b^{\re}_+({R})-c^{\re}({R}),
\end{equation}
which is the Klein formula in its original form \cite{Kl-paper}.

\subsection{Example}
Let us go back to the log-front from Figure \ref{froz}.
By Theorem \ref{Harthm} the sum of the real cusps of the
log-front and twice the number of real solitary nodes equals 24.
Indeed, we have $n=6$, $\Vol(\Delta_P,\Delta_Q)=8$,
$\Area(\Delta_P)=4$, $|\dd\Delta_Q|=12$, $\chi(\tilde{Q})=2$,
$2b_+^{\re}(Q)=b_+(\dd Q)=c(Q)=0$ while $\epsilon(\Delta_P,\Delta_Q)=48$.

\begin{figure}[!h]
  \centering
  {\scalebox{0.5}{\includegraphics{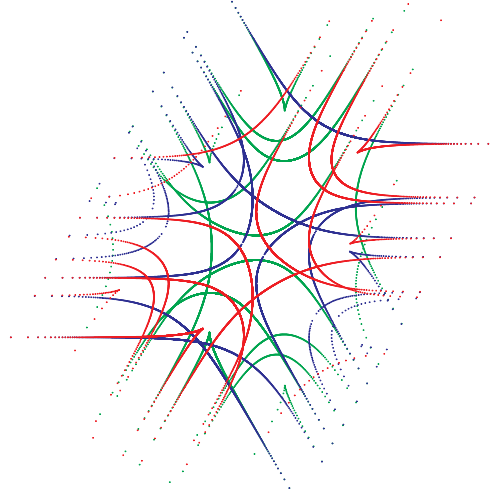}}}
  \caption{Images in the other three quadrants}
  \label{3quad}
\end{figure}

There are 6 real cusps visible on Figure \ref{froz}. Furthermore,
Figure \ref{3quad} shows images of the remaining 3 quadrants
under the map $(x,y)\mapsto (\log|x|,\log|y|)$. There are 12 more
real cusps. Thus by Theorem \ref{Harthm} our log-front has 3 real
solitary points.

%%%%%%%%%%%%%%%%%%%%%%%%%%%%%%%%%%%%%%%%%%%%%%%%%%%%%%%%%%%%%%%%%%%%%%%%%%%%

\newpage
\section*{Appendix A: Harnack curves and their algae.}

In his 1876 paper \cite{Ha} A. Harnack produced
for each $d$ examples of algebraic curves of degree $d$ in $\rp^2$
with $\frac{(d-1)(d-2)}{2}+1$ real (topological) components.
Furthermore, in the same paper he has shown that $\frac{(d-1)(d-2)}{2}+1$
is the upper bound for the number of components
of any curve of degree $d$ in $\rp^2$

Among the curves constructed by Harnack there were
some ``canonical" curves whose topological arrangement
in $\rp^2$ is especially easy to describe.
Note that any component of a smooth curve is either
contractible (i.e. bounds a disk in $\rp^2$)
or is isotopic to $\rp^1\subset\rp^2$. A contractible component
is called an {\em oval} while the disk bounded by it is called
the {\em interior} of the oval. The oval whose interior is
disjoint from other ovals is called {\em empty}.

For an odd $d=2k+1$
there is a smooth algebraic curve of degree $d$ that
consists of $k(2k-1)$ empty ovals
and a non-contractible component.
For an even $d=2k$ there is a smooth algebraic curve of degree $d$ that
consists of $(k-1)(2k-1)$ empty ovals and
one other oval whose interior contains $\frac{(k-1)(k-2)}{2}$
of the $(k-1)(2k-1)$ empty ovals.

%Let $\R A\subset\tor$ be a curve whose Newton polygon is
%a triangle with vertices $(0,0)$, $(d,0)$ and $(0,d)$.
Let $P\subset\cT$ be a curve given by a real polynomial
with the Newton Polygon $\Delta_P$ and $\cT$ be the toric
surface corresponding to $\Delta_P$.
Since $P$ is defined over
$\R$ there is the real locus $\R P=P\cap\rtor$ which coincides
with the fixed point set of the involution of complex conjugation
on $P$.
%Consider the logarithm map $$\Log:\tor\to\R^2,$$
%$\Log(z,w)=(\log|z|,\log|w|)$. The image $\Log(A)\subset\R^2$
%is called the amoeba of $A$, see \cite{GKZ}.

\begin{defn}\label{sH}
A curve $P$ is called a {\em Harnack
curve} if for every $(x,y)\in\R^2$ the set $\Log^{-1}(x,y)\cap P$
consists of no more than two points.
%there exist three arcs $\alpha_1,\alpha_2,\alpha_3$ on the
%same topological component of $\R A$ such that
\end{defn}

\begin{rmk}
Earlier we called such curves {\em simple Harnack curves}
to distinguish them from other curves in the Harnack construction.
However, by now we have convinced ourselves that these {\em simple} curves
are the most beautiful in the Harnack series of constructions.
We propose to drop ``simple" from their name and call them
{\em Harnack curves}.
\end{rmk}

Harnack curves (from Definition \ref{sH})
exist for any convex lattice polygon $\Delta$, cf. \cite{IV}.
In \cite{Mi} it was shown that the topological type of
the triad $(\R\cT;\R P,\dd\R\cT)$
%(and of the compactification $(\R\c T,\R {A})$
%in the real toric surface $\R T_{\Delta}$ corresponding to $\Delta$)
depends only on $\Delta$ if $\R P$ is a smooth Harnack curve
transverse to infinity.
%Furthermore, as usual we may compactify the torus $\rtor$ to the toric surface
%$\R\cT$ corresponding to the Newton polygon $\Delta$ and the topological
%type of $(\R\c T,\overline{\R P})$ will remain determined only by $\Delta$.
%As in the main part of the paper from now on in the Appendix
%we denote with $P$ the corresponding compactification of the curve
%so that $P\setminus\dd P=P\cap\tor$.

In the case when $\Delta$ is
a triangle with vertices $(0,0)$, $(d,0)$ and $(0,d)$ we have
$\R\cT=\rp^2$ and
%the compactification $\R\bar{A}$ of
%$\R A$ in $\rp^2$
$\R P\subset\rp^2$ is a curve consisting of $\frac{(d-1)(d-2)}{2}$
empty oval and one other component which is non-contractible if
$d$ is odd and an oval containing $\frac{(d-2)(d-4)}{8}$ empty ovals
if $d$ is even. In \cite{KO1} it was shown that all such curve
form a contractible subspace in the space of all real curves
of degree $d$.

Recently it was discovered that
Harnack curves possess many extremal characteristic properties, among them
are the following.
\begin{itemize}

\item (\cite{MR}) We have
$$\Area(\Log(P))=\pi^2\Area(\Delta),$$
if $\R P$ is a Harnack curve with the Newton polygon $\Delta$.
In the same time by \cite{PR}
for any curve $A\subset\tor$ given by a (not necessarily
real) polynomial with the Newton polygon $\Delta$ we have
$$\Area(\Log(A))\le\pi^2\Area(\Delta).$$

Furthermore, if
$\Area(\Log(A))=\pi^2\Area(\Delta)$ then $A$ can be translated
(by a multiplication  with some $(z_0,w_0)\in\tor$) to a
Harnack curve, see \cite{MR}.

\item (\cite{Mi}) The curve $\Log(\R P)\subset\R^2$ is embedded,
does not have inflection points and contains $\#(\Int(\Delta)\cap\Z^2)$
compact components (called {\em ovals}) if $\R P$ is a Harnack curve.
Each such oval comes from one of the four quadrants in $\rtor$.
The number of the ovals coming from the four quadrant
equals to the number of lattice points $(j,k)\in\Z^2$ in the interior
$\Int(\Delta)$ with given residue mod 2: $j\equiv j_0\pmod{2}$,
$k\equiv k_0\pmod{2}$ (clearly there are four possible pairs of residues).

Furthermore any curve with that many ovals of $\Log(\R P)\approx P\setminus\dd P$
is Harnack if the remaining component of $\R P$ intersects the infinity $\dd\R\cT$
in a maximal way (see \cite{Mi}).
\end{itemize}

\begin{figure}[h]
\centerline{\psfig{figure=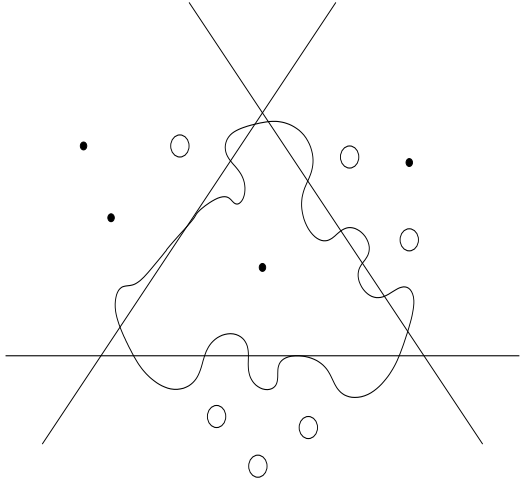,height=2.5in,width=2.7in}}
\caption{\label{exaHarnack} A singular simple Harnack curve in $\rp^2$ and
its position with respect to the boundary divisor.}
\end{figure}

The only singularities of a Harnack curve
$\R P$ are isolated double points
(the singularities of type $A_1^+$ according to \cite{AVGZ}) in $\rtor$,
see \cite{MR}.
Note that even though $\R P$ has to be smooth near the boundary
divisor $\dd\R\cT$ it does not have to be transverse
to the boundary divisor. An example of a singular simple Harnack
curve of degree 6 in $\rp^2$ which is not transverse to the boundary
divisor is sketched in Figure \ref{exaHarnack}.
%\end{itemize}

The goal of this appendix is to give yet another characteristic property
of the Harnack curves in terms of their algae.
Suppose that $P$ is a Harnack curve %, see \cite{Mi},\cite{MR},\cite{KO}.
Recall that in the previous section we denoted with $\hat{P}$
the result of the real blowing up of $\tilde{P}$ at the points $\dd P$
and with $$\Alga_P:\hat{P}\to T$$
such that $\Alga_P|_{\tilde{P}\setminus\dd P}=
\beta\circ\Arg|_{\tilde{P}\setminus\dd P}$.
%of $\tilde{P}\setminus\nu_{\bar{P}}^{-1}(\bar{P}\setminus P)$
%and $\R\hat{P}$ is the real part of $\hat{P}$.
%Consider the argument map $\Arg:\tor\to T=S^1\times S^1$,
%$\Arg(z,w)=(\arg(z),\arg(w)$.
%Clearly, $\Arg|_P$ extends
%to the map $$\hat\Arg_P:\hat{P}\to T.$$
Recall that $P_0=\Alga^{-1}_P(0_T)$.

\begin{lemmaA}\label{arg-harnack}
If $P$ is a Harnack curve then the restriction of the map $\Alga_P$
to $\hat{P}\setminus P_0$
$$\hat{P}\setminus P_0\to T\setminus\{0_T\}$$ is
an unbranched covering of degree $2Area(\Delta_P)$.
%Furthermore, $P_0\setminus\R\hat{P}$ is a finite set.
\end{lemmaA}
%Recall (see section \ref{logGauss}) that $a_P$ is the degree of
%the logarithmic Gauss map of $P$. If $\bar{P}$ is a smooth curve
%orthogonal to the boundary divisor $\C T_{\Delta_P}\setminus\tor$
%then $a_P=\Area\Delta_P$.
\begin{proof}
%[Proof of Lemma A1]
The critical points of the map $\Arg|_P$ are the points
such that the image of the Logarithmic Gauss map $\gamma_P$ is in $\rp^1$
(see \cite{Mi}).
By \cite{Mi} $\gamma^{-1}(\rp^1)=\R\tilde{P}$.
The map $\Alga_P|_{\hat{P}\setminus P_0}$ is proper since
$\Alga_P^{-1}(0_T)= P_0$.

Since $\chi(T\setminus\{0_T\})=-1$ we have the degree of the covering
$\Alga_P|_{\hat{P}\setminus P_0}$ equal to $-\chi(P\setminus P_0)$.
Let us start with a generic $P$ and
study how does $\chi(P\setminus P_0)$ changes when we
deform $P$ in the class of Harnack curve.
Since $P$ is a simple Harnack curve we have $P\setminus P_0$ non-singular.
The only singularities of $P$ are the real isolated double points.
Each such point contributes $2$ to the Euler characteristic of the
normalization of $P$, but also $2$ gets subtracted when we remove
this point from $\hat{P}$.

Note that $P$ may also have the ``boundary" singularity.
This means that ${P}$ is not transversal to the boundary divisor.
In this case each point of tangency of order $m$ between ${P}$
and $\dd\cT$ gives $m-1$ points in $\hat{P}\setminus P_0$.
Thus, while $\chi(\hat{P})$ gets increased by $m-1$
in the case of such tangency (in comparison with $\chi(\hat{P})$
in the transversal case) in turn $\chi(\hat{P}\setminus P_0)$ gets
decreased by $m-1$. Thus in both cases the net effect
of possible singularities on $\chi(\hat{P}\setminus P_0)$ is zero,
so for the computation of the degree of our covering
we may assume that $\tilde{P}$ is smooth and transversal to
the boundary divisor of $\cT$.
In this case $$\chi(\hat{P}\setminus P_0)=\chi(\hat{P})=\chi(P)=\Area\Delta.$$
The last equality is a corollary of
Khovanskii's formula \cite{Kh}.
%To compute the degree of the covering $\Arg_P$ we note that
%%by Khovanskii's formula \cite{Kh} the Euler characteristic
%%of $\hat{P}$ (that is equal to that of $P$) is $\Area(\Delta_P)$.
%$\chi(\hat{P})$ is equal to $a_P$ since the logarithmic Gauss map
%has no real critical points (see \cite{Mi}).
\end{proof}

We may compactify the set-up of Lemma \ref{arg-harnack}
to get the following Theorem describing the alga of
a Harnack curve.
Denote the blow-up of $T$ centered in $0_T$ with
$B:\hat{T}\to T$.

Note that $\hat{T}$ is a surface with a natural involution
induced by $\conj$: $(\alpha,\beta)\mapsto (-\alpha,-\beta)$,
here we think of $\alpha,\beta\in S^1$ as arguments of complex numbers.
Clearly, the 1-dimensional part of the
fixed-point set of this involution is the exceptional
divisor of the blow-up. Denote it with $\R\hat{T}$.

For the compactifying theorem we need to blow up $\hat{P}$
even further. Recall that $\hat{P}$ is a smooth (real) surface
equipped with an involution coming from complex conjugation.
Note that $0_T$
is the regular value of a map $P\setminus\R P\to T$
and thus $P_0\subset\hat{P}$
is a disjoint union of embedded circles corresponding to the
ovals of $\R\tilde{P}$ and some isolated points.
Let $\Pi$ be the result of (real) blow-up of $\hat{P}$
at $P_0$. Clearly, a blow-up at a smooth submanifold of codimension 1
does not change the surface $\hat{P}$ thus only blowups at
the isolated points of $P_0$ matter.
Such points come either from isolated double points of $P$
or from tangency of $\tilde{P}$ with $\dd\cT$.
%the boundary divisor of $\C T_{\Delta}$.
Note that
$$\chi(\Pi)=\chi(\hat{P}\setminus P_0).$$

\begin{theoremA}
There exists a map $p$ completing the commutative diagram
\begin{center}\mbox{}
\xymatrix{
\Pi \ar[r]^{p} \ar[dr]_{\Alga_P} & \hat{T}\ar[d]^{B}\\
  & T}.
\end{center}
The map $p$ is a covering of degree $2\Area(\Delta_P)$
which is equivariant with respect to the involution of complex conjugation
defined on $\Pi$.
%In particular, $$\R\hat{P}=\hat\Arg_P^{-1}(\R\hat{T}).$$

Furthermore, any real curve $P$ such that $\Alga_P$ lifts to
a covering $\Pi\to\hat{T}$ is a Harnack curve.
\end{theoremA}
\begin{proof}
The point of $\R\hat{T}$ is specified by a tangent line to $T$ at $0_T$.
The logarithmic Gauss map $\gamma_P$ takes real values
at $\R\hat{P}$. Thus the tangent line at $s\in\R\hat{P}$
is real and gives (after multiplication by $i$)
a tangent direction at $0_T\in T$. We define the value $p(s)$
to be this direction. At points of $P_0\setminus \R\hat{P}$
the map is defined by the blowup itself (recall that $0_T$
is the regular value of the map $P\setminus\R P\to T$).

For the converse we note that by \cite{Mi} the amoeba map
$\Log|_P$ does not have critical points outside of $\R P$
since $\Alga_P$ does not have any. Thus $\gamma_P$ cannot
have any critical points on $\R P$.

\end{proof}

\end{document}